\documentclass[12pt]{amsart}
\usepackage{latexsym}
\usepackage{amsmath,amsfonts,amssymb}
\numberwithin{equation}{section}
\newtheorem{theo}{Theorem}[section]

\newtheorem{prop}[theo]{Proposition}
\newtheorem{rem}[theo]{Remark}

\newcommand{\real}{\mathbb R}
\newcommand{\Nat}{\mathbb N}

\newcommand{\diam}{\operatorname{diam}}

\newcommand{\bequ}{\begin{equation}} %
\newcommand{\eequ}{\end{equation}} %
\newcommand{\bequs}{\begin{equation*}} %
\newcommand{\eequs}{\end{equation*}} %
\newcommand{\vphi}{\varphi}
\newcommand{\rondu}{\mathcal U}%
\begin{document} %

\title[Compact operators]{ Compact operators on spaces with asymmetric norm}

\author{  S. Cobza\c s}

\address{\it Babe\c s-Bolyai University, Faculty of Mathematics
and Computer Science, Ro-3400 Cluj-Napoca, Romania,\;\; E-mail:
scobzas@math.ubbcluj.ro} %

\begin{abstract}

The aim of the present paper is to define compact operators on
asymmetric normed spaces and to study some of their properties.
The dual of a bounded linear operator is  defined and a Schauder
type theorem is proved within this framework. The paper contains
also a short discussion on various completeness notions for
quasi-metric and for quasi-uniform spaces.

 AMS 2000 MSC: Primary: 47B07;\;
 Secondary: 46B99; 46S99; 47A05; 54E15; 54E25

\end{abstract}
\maketitle
\section{Introduction}

An asymmetric norm on a real vector space $X$ is a functional
$p:X\to [0,\infty)$ satisfying the conditions %
$$ %
\mbox{(AN1)}\; p(x)=p(-x)=0\Rightarrow x=0;\quad\mbox{(AN2)}\;
p(\alpha x)=\alpha p(x);\quad \mbox{(AN3)}\; p(x+y)\leq p(x)+p(y), %
$$%
for all $x,y\in X$ and $\alpha \geq 0.$ A quasi-metric on a set
$X$ is a mapping $\rho:X\times X\to [0,\infty)$ satisfying the
conditions %
$$%
\mbox{(QM1)}\; \rho(x,y)=\rho(y,x)=0 \iff x=y ;\quad\mbox{(QM2)}\;
\rho(x,z)\leq \rho(x,y)+\rho(y,z), %
$$%
for all $x,y,z \in X.$  If the mapping $\rho$ satisfies only the
conditions $\,\rho(x,x)=0,\, x\in X,\,$ and  (QM2), then it is
called a {\it quasi-pseudometric}. If $p$ is an asymmetric norm on
a vector space $X,$ then the pair $(X,p)$ is called an asymmetric
normed space. Similarly, $(X,\rho)$ is called a quasi-metric
space. If $p$ is an asymmetric norm on a vector space $X$, then
$\rho(x,y)=p(y-x),\;x,y\in X,$ is a quasi-metric on $X.$ A closed,
respectively open, ball in a quasi-metric space is defined
by %
\bequs %
B_\rho(x,r)=\{y\in X : \rho(x,y)\leq r\},\quad
B'_\rho(x,r)=\{y\in X : \rho(x,y)< r\}, %
\eequs %
for $x\in X$ and $r>0.$ In the case of an asymmetric norm $p,\,$
one denotes by $B_p(x,r), B'_p(x,r)$ the corresponding balls and
by $B_p = B_p(0,1), B'_p=B'_p(0,1), $ the unit balls. In this case
the following equalities hold %
\bequs %
B_p(x,r)=x+rB_p \quad\mbox{and}\quad B_p'(x,r)=x+rB'_p. %
\eequs %

The family of sets $\,B'_\rho(x,r),\, r>0, $  is a base of
neighborhoods of the point $x\in X$ for the topology $\tau_\rho$
on $X$ generated by the quasi-metric $\rho.$ The family
$\,B_\rho(x,r),\, r>0, $ of closed balls is also a neighborhood
base at $x$ for $\tau_\rho.$

A quasi-uniformity on a set $X$ is a filter $\mathcal U$ such that %
\bequs %
\begin{aligned} %
\mbox{(QU1)}&\qquad \Delta(X)\subset U,\; \forall U\in \rondu;\\
\mbox{(QU1)}&\qquad \forall U\in \rondu,\; \exists V\in \rondu,\;
\mbox{such that } \; V\circ
V\subset U, %
\end{aligned}
\eequs %
where $\Delta(X)=\{(x,x) : x\in X\}$ denotes the diagonal of $X$
and, for $M,N\subset X\times X,$  %
\bequs %
M\circ N =\{(x,z)\in X\times X : \exists y\in X,\; (x,y)\in
M\;\mbox{and} \; (y,z)\in N\}. %
\eequs %

If the filter $\rondu$ satisfies also the condition %
\bequs %
\mbox{(U3)}\qquad \forall U,\; U\in \rondu\;\Rightarrow\;
U^{-1}\in\rondu,
\eequs %
where %
$$%
U^{-1}=\{(y,x)\in X\times X : (x,y)\in U\}, %
 $$ %
then $\rondu$ is called a {\it uniformity} on $X.$  The sets in
$\rondu$ are called {\it entourages} (or {\it vicinities}).

For $U\in \rondu, \, x\in X$ and $Z\subset X$ put %
$$%
U(x) =\{y\in X: (x,y)\in U\}\quad\mbox{and}\quad U[Z]=\cup\{U(z):
z\in Z\}. %
$$ %
 A quasi-uniformity $\rondu$ generates a topology $\tau(\rondu)$
 on $X$ for which the family of sets%
 \bequs %
\{ U(x): U\in \rondu\} %
\eequs %
is a base of neighborhoods of the point $x\in X.$  A mapping $f$
between two quasi-uniform spaces $(X,\rondu),\,(Y,\mathcal W)$ is
called {\it quasi-uniformly continuous} if for every $W\in
\mathcal W$ there exists $U\in \rondu$ such that $(f(x),f(y))\in
W$ for all $(x,y)\in U.$ By the definition of the topology
generated by a quasi-uniformity, it is clear that a
quasi-uniformly continuous mapping is continuous with respect to
the topologies $\tau(\rondu),\, \tau(\mathcal W).$

If $(X,\rho)$ is a quasi-metric space, then %
\bequs %
 B'_\epsilon =\{(x,y)\in X\times X : \rho(x,y)< \epsilon\},\;
 \epsilon > 0, %
 \eequs %
 is a basis for a quasi-uniformity $\rondu_\rho$ on $X.$  The
 family %
 \bequs %
 B_\epsilon =\{(x,y)\in X\times X : \rho(x,y)\leq \epsilon\},\; \epsilon > 0, %
 \eequs %
 generates the same quasi-uniformity. The topologies generated by the
 quasi-metric $\rho$ and by the quasi-uniformity $\rondu_\rho$
 agree, i.e., $\, \tau_\rho=\tau(\rondu_\rho).$

 The lack of the symmetry, i.e., the omission of the axiom (U3),
makes the theory of quasi-uniform spaces to differ  drastically
from that of uniform spaces. An account of the theory up to 1982
is given in the book by Fletcher and Lindgren \cite{FL}. The
survey papers by K\" unzi \cite{kunz93a,kunz93b,kunz01,kunz02b}
are good guides for subsequent developments. Another book on
quasi-uniform spaces is \cite{MN}.

On the other hand, the theory of asymmetric normed spaces has been
developed in a series of papers
\cite{cob04,cob-mus04,rom-raf-per02a,rom-raf-per02b,rom-raf-per03a,rom-raf-per03b,rom-raf-per04b},
following ideas from the theory of (symmetric) normed spaces and
emphasizing similarities as well as differences between the
symmetric and the asymmetric case.

Let $(X,p)$ be an asymmetric normed space. The functional $\bar
p(x)=p(-x),\, x\in X,$ is also an asymmetric norm on $X,$ called
the conjugate of $p,\; p_s(x)=\max\{p(x),\bar p(x)\},\; x\in X,$
is a (symmetric)
norm on $X$ and the following inequalities hold %
\bequs %
|p(x)-p(y)|\leq p_s(x-y) \quad\mbox{and}\quad|\bar p(x)-\bar
p(y)|\leq
p_s(x-y),\; \forall x,y\in X.%
\eequs %

For a quasi-metric space one defines similarly the conjugate of
$\rho$ by $\bar \rho(x,y)=\rho(y,x)$ and the associated
(symmetric) metric by $\rho_s(x,y)=\max\{\rho(x,y),\rho(y,x)\}, $
for $x,y\in X.$

Let $(X,p),\, (Y,q)$ be two asymmetric normed space. A linear
mapping $A:X\to Y$ is called {\it bounded},  ($(p,q)$-bounded if
more precision is needed), or {\it semi-Lipschitz}, if there
exists a
number $\beta \geq 0$ such that %
\bequ\label{def.sL} %
q(Ax)\leq \beta p(x), %
\eequ %
for all $x\in X.$ The number $\beta$ is called  a semi-Lipschitz
constant for $A.$ For properties of semi-Lipschitz functions and
of spaces of semi-Lipschitz functions see
\cite{mus01,mus02,rom-sanc00,rom-sanc05}.

 The operator $A$ is continuous with respect
to the topologies $\tau_p,\tau_q$ ($(\tau_p,\tau_q)$-continuous)
if and only if it is bounded and if and only if it is
quasi-uniformly continuous with respect to the quasi-uniformities
$\rondu_p$ and $\rondu_q$ (see \cite{aleg-fer-greg93} and
\cite{rom-raf-per03a}). Denote by $(X,Y)^\flat_{p,q},$ or simply
by $(X,Y)^\flat$ when there is no danger of confusion, the set of
all $(p,q)$-bounded linear operators.  The set $(X,Y)^\flat$ need
not be a linear subspace but merely a convex cone in the space
$(X,Y)^\#$ of all linear operators from $X$ to $Y,$ i.e., $A+B\in
(X,Y)^\flat$ and $\alpha A\in (X,Y)^\flat, $ for any $A,B\in
(X,Y)^\flat$ and $\alpha \geq 0.$ Following \cite{rom-raf-per03a},
we shall call
$(X,Y)^\flat$ a {\it semilinear space}.  The functional %
\bequ\label{def.as-n} %
\|A|=\|A|_{p,q}=\sup\{q(Ax) : x\in B_p\} %
\eequ %
is an asymmetric norm on the semilinear space $(X,Y)^\flat,$ and
$\|A|$ is the smallest semi-Lipschitz constant for $A,$ i.e., the
smallest number for which the inequality \eqref{def.sL} holds.

Denote by $(X,Y)^*_s$ the space of all continuous linear operators
from $(X,p_s)$ to $(Y,q_s),$ normed by %
\bequ\label{def.nA} %
\|A\|=\|A\|_{p_s,q_s}=\sup\{q_s(Ax):x\in X, \; p_s(x)\leq 1\}, \;\; A\in (X,Y)^*_s. %
\eequ %

It was shown in \cite{rom-raf-per03a} that
$(X,Y)_{p,q}^\flat\subset (X,Y)^*_s, $ and $\|A|\leq \|A\|$ for
any $A\in (X,Y)^\flat_{p,q}.$

Consider on $\real$ the asymmetric norm $u(\alpha) =
\max\{\alpha,0\}, \, \alpha\in \real.$ Its conjugate is $\bar
u(\alpha)=\max\{-\alpha,0\}$ and $u_s(\alpha)=|\alpha|$ is the
absolute value norm on $\real.$ The topology $\tau_u$ on $\real$
generated by $u,$  called the upper topology of $\real$, has as
neighborhood basis of a point $\alpha\in \real$ the family of
intervals $(-\infty,\alpha+\epsilon),\,\epsilon >0.$

The space of all linear bounded functionals from an asymmetric
normed space $(X,p)$ to $(\real,u)$ is denoted by $X^\flat_p\,.$
Notice that, due to the fact that $p$ is non-negative, we have
$$%
\forall x\in X, \; u(\vphi(x))\leq \beta p(x) \;\iff\; \; \vphi(x)\leq \beta p(x), %
$$ %
for any linear functional  $\vphi:X\to\real,$ so the asymmetric
norm of a
functional $\vphi\in X^\flat_p$ is given by %
$$%
\|\vphi|=\|\vphi|_p=\sup\{\vphi(x) : x\in X,\; p(x)\leq 1\}. %
$$%

Also, the continuity of $\vphi$  from $(X,\tau_p)$ to
$(\real,\tau_u)$ is equivalent to its upper semi-continuity from
$(X,\tau_p)$ to $(\real,|\;|),$ (see
\cite{aleg-fer-greg97,aleg-fer-greg99,aleg-fer-greg93}).

In \cite{rom-raf-per03a} it was defined the analog of the
$w^*$-topology on the space $X^\flat_p,$  which we denote by
$w^\flat,\,$ having as a base of $w^\flat$-neighborhoods of an
element $\vphi_0\in
X^\flat_p$  the sets %
\bequ\label{def.w-flat} %
V_{x_1,...,x_n;\,\epsilon}(\vphi_0)=\{\vphi\in X^\flat_p :
\vphi(x_i)-\vphi_0(x_i)\leq \epsilon,\; i=1,...,n\},  %
\eequ%
for $n\in \Nat,\, x_1,...,x_n\in X,\,$ and $\epsilon >0.$

Since %
$$%
V_{x;\,\epsilon}(\vphi_0)\cap
V_{-x;\,\epsilon}(\vphi_0)=\{\vphi\in X^\flat_p :
|\vphi(x)-\vphi_0(x)|\leq \epsilon\}, %
$$%
it follows that the topology $w^\flat$ is the restriction to
$X^\flat$ of the $w^*$-topology of $X^*_s =(X,p_s)^*.$

Some results on $w^\flat$-topology were proved in
\cite{rom-raf-per03a} as, for instance, the analog of the
Alaoglu-Bourbaki theorem: the
polar %
\bequ \label{def.polar} %
B^\flat_p=\{\vphi \in X^\flat : \vphi(x)\leq 1, \; \forall x\in
B_p \}%
\eequ %
of the unit ball $B_p$  of $(X,p)$ is $w^\flat$-compact. Other
results on asymmetric normed spaces, including separation of
convex sets by closed hyperplanes  and a Krein-Milman type
theorem, were obtained in \cite{cob04}. Asymmetric locally convex
spaces were considered in \cite{cob05}.  Best approximation
problems in asymmetric normed spaces were studied in \cite{cob04}
and \cite{cob-mus04}.

 The topology $w^\flat$ is derived from a quasi-uniformity
 $\mathcal W^\flat_p$ on $X^\flat_p$ with a basis formed of the sets%
 \bequ\label{def.w-flat-unif} %
 V_{x_1,...,x_n;\,\epsilon}=\{(\vphi_1,\vphi_2)\in X^\flat_p \times X^\flat_p :
\vphi_2(x_i)-\vphi_1(x_i)\leq \epsilon,\; i=1,...,n\},  %
\eequ%
for $n\in \Nat,\, x_1,...,x_n\in X$ and $\epsilon >0.$ Note that,
for fixed $\vphi_1=\vphi_0,$ one obtains the neighborhoods from
\eqref{def.w-flat}.

On the space $(X,Y)^*_s$ we shall consider several
quasi-uniformities. Namely, for $\mu\in \{p,\bar p,p_s\}$ and
$\nu\in \{q,\bar q,q_s\}$ let $\rondu_{\mu,\nu}$ be the
quasi-uniformity
generated by the basis %
\bequ\label{def.q-u} %
U_{\mu,\nu;\,\epsilon} =\{(A,B) ; A,B\in (X,Y)^*_s, \;
\nu(Bx-Ax)\leq
\epsilon,\;\forall x\in B_\mu,\}, \;\;\epsilon >0,%
\eequ %
where $B_\mu=\{x\in X : \mu(x)\leq 1\}$ denotes the unit ball of
$(X,\mu).$ The induced quasi-uniformity on the semilinear subspace
$(X,Y)^\flat_{\mu,\nu}$ of $(X,Y)^*_s$ is denoted also by
$\rondu_{\mu,\nu}$ and the corresponding topologies  by
$\tau(\mu,\nu).$ The uniformity $\rondu_{p_s,q_s}$ and the
topology $\tau(p_s,q_s)$ are those  corresponding to the norm
\eqref{def.nA} on the space $(X,Y)^*_s.$

In the case of the dual space $X^\flat_\mu$ we shall use the
notation $\rondu^\flat_\mu$ for the quasi-uniformity $\rondu_{\mu,u}\,.$ %

\section{Completeness and compactness in quasi-metric and in
quasi-uniform spaces}

The lack of symmetry in the definition of quasi-metric and
quasi-uniform spaces causes a lot of troubles, mainly concerning
completeness, compactness and total boundedness in such spaces.
There are a lot of completeness notions in quasi-metric and in
quasi-uniform spaces, all agreeing with the usual notion of
completeness in the case of metric or uniform spaces,  each of
them having its advantages and weaknesses.

We shall describe briefly some of these notions along with some of
their properties.

The first one is that of bicompleteness. A quasi-metric space
$(X,\rho)$ is called {\it bicomplete} if the associated symmetric
metric space $(X,\rho_s)$ is complete. A bicomplete asymmetric
normed space $(X,p)$ is called also a {\it biBanach space.} The
existence of a bicompletion of an asymmetric normed space was
proved in \cite{rom-raf-per02a}. The notion can be considered also
for an {\it extended} (i.e. taking values in $[0,\infty]$)
quasi-metric, or for an extended asymmetric norm on a semilinear
space.

In \cite{rom-raf-per03a} it was defined an extended asymmetric
norm on $(X,Y)^*_s$ by %
\bequ\label{def.anA} %
\|A|^*_{p,q} =\sup\{q(Ax) :x\in B_p\}, \;\; A\in (X,Y)^*_s. %
\eequ %

The identity mapping $\operatorname{id}_\real$ is continuous from
$(\real,u)$ to $(\real,u),$ but for $-\operatorname{id}_\real$ we have %
$$%
\|-\operatorname{id}_\real|^*_{u,u}=\sup\{-\alpha : u(\alpha)\leq
1\}\geq\sup\{-\alpha:\alpha \leq 0\} =+\infty, %
$$%
because $u(\alpha)=0\leq 1$ for $\alpha \leq 0. $ It follows that
$\|A|^*_{p,q}$ can take effectively the value $+\infty.$

 If the asymmetric
normed space $(Y,p)$ is bicomplete, then  the space $(X,Y)^*_s$ is
complete with respect to the symmetric extended norm
$(\|\,|^*_{p,q})_s$ and $(X,Y)^\flat_{p,q}$ is a
$(\|\,|^*_{p,q})_s$-closed semilinear subspace  of $(X,Y)^*_s,$ so
it is $\|\,|_{p,q}$-bicomplete (see \cite{rom-raf-per03a}). %

In the case of a quasi-metric space $(X,\rho)$ there are also
other completeness notions. We present them following
\cite{reil-subram82}, starting with the definitions of Cauchy
sequences.

A sequence $(x_n)$ in $(X,\rho)$ is called

(a) \;{\it left (right) $\rho$-Cauchy} if for every $\epsilon > 0$
there exist $x\in X$ and $n_0\in \Nat$ such that

\quad\quad $\forall n\geq n_0,\; \rho(x,x_n)<\epsilon$
(respectively $\rho(x_n,x)<\epsilon$)\,;

(b) \; $\rho$-{\it Cauchy} if for every $\epsilon >0$ there exists
$n_0\in \Nat$ such that

\quad\quad $\forall n,k\geq n_0,\; \rho(x_n,x_k)< \epsilon\,$;

(c)\;  {\it left (right)-K-Cauchy} if  for every $\epsilon >0$
there exists $n_0\in \Nat$ such that

\quad \quad$\forall n,k, \; n\geq k\geq n_0\; \Rightarrow \;
\rho(x_k,x_n)< \epsilon$ (respectively
$\rho(x_n,x_k)<\epsilon$)\,;

(d)\; {\it weakly left(right) K-Cauchy} if for every $\epsilon >0$
there exists $n_0\in \Nat$ such that

\quad \quad$\forall n\geq n_0,\; \rho(x_{n_0},x_n)< \epsilon$
(respectively $\rho(x_n,x_{n_0})< \epsilon$).

These notions are related in the following way: \vspace*{1em}

left(right) $K$-Cauchy $\;\Rightarrow\; $ weakly left(right)
$K$-Cauchy $\;\Rightarrow\;$ left(right) $\rho$-Cauchy,
\vspace*{1em}\newline%
and no one of the above implications is reversible
(see \cite{reil-subram82}).

Furthermore, each $\rho$-convergent sequence is $\rho$-Cauchy, but
for each of the other notions there are examples of
$\rho$-convergent sequences that are not Cauchy, which is a major
inconvenience. Another one is that closed subspaces of complete
(in some sense) quasi-metric spaces need not be complete. If each
convergent sequence in a regular quasi-metric space $(X,\rho)$
admits a left $K$-Cauchy subsequence, then $X$ is metrizable
(\cite{kunz-reil93}. This result shows that putting too many
conditions on a quasi-metric, or on a quasi-uniform space, in
order to obtain results similar to those in the symmetric case,
there is the danger to force the quasi-metric to be a metric and
the quasi-uniformity a uniformity. In fact, this is a general
problem when dealing with generalizations.

For  each of these notions of Cauchy sequence one obtains a notion
of sequential  completeness,  by asking that each corresponding
Cauchy sequence be convergent in $(X,\rho).$ These notions of
completeness are related in the
following way: %
\vspace*{1em} %

left (right) $\rho$-sequentially complete $\;\Rightarrow\;$ weakly
left (right) $K$-sequentially complete $\;\Rightarrow\;$

$\;\Rightarrow\;\rho$-sequentially complete.

\vspace*{1em}

In spite of the obvious fact that left $\rho$-Cauchy is equivalent
to right $\bar \rho$-Cauchy, left $\rho$- and right
$\bar\rho$-completeness do not agree, due to the fact that right
$\bar \rho$-completeness means that every left $\rho$-Cauchy
sequence converges in $(X,\bar \rho),$ while left
$\rho$-completeness means the convergence of such sequences in the
space $(X,\rho).$ For concrete examples and counterexamples, see
\cite{reil-subram82}.

A subset $Y$ of a quasi-metric space $(X,\rho)$ is called {\it
precompact} if for every $\epsilon > 0$ there exists a finite
subset $Z$ of $X$ such that %
$$%
Y\subset \cup\{B_\rho(z,\epsilon) : z\in Z\}. %
$$%

The set $Y$ is called {\it totally bounded } if for every
$\epsilon> 0,\; Y$ can be covered by a finite family of sets of
diameter less that $\epsilon,$ where the diameter of a subset $A$
of $X$ is defined by %
$$%
\diam (A) =\sup\{\rho(x,y) : x,y\in A\}. %
$$%

As it is known, in metric spaces the precompactness and the total
boundedness are equivalent notions, a result that is not longer
true in quasi-metric spaces, where  precompactness is strictly
weaker than total boundedness, see \cite{lambrin77} or \cite{MN}.

In spite of these peculiarities there are some positive results
concerning Baire theorem and compactness. For instance, any
compact quasi-metric space is left $K$-sequentially complete and
precompact.  If $(X,\rho)$ is precompact and left
$\rho$-sequentially complete, then it is sequentially compact (see
\cite{fer-greg85,reil-subram82}). Hicks \cite{hiks88} proved some
fixed point theorems in quasi-metric spaces (see also
\cite{hiks71,hiks80})

K\" unzi et al \cite{kunz-reil93} proved that a quasi-metric space
is compact if and only if it is precompact and left
$K$-sequentially complete, and studied the relations between
completeness, compactness, precompactness, total boundedness  and
other related notions in quasi-uniform spaces.

Notice also that in quasi-metric spaces compactness, countable
compactness  and sequential compactness are different notions (see
\cite{fer-greg83} and \cite{kunz83}).

The considered completeness notions can be extended to
quasi-uniform spaces by replacing sequences by filters or nets
(for nets, see \cite{sunderh95,sunderh97}). Let $(X,\rondu)$ be a
quasi-uniform space, $\rondu^{-1}=\{U^{-1} : U\in \rondu\}$ the
conjugate quasi-uniformity on $X,$ and $\rondu_s=\rondu \vee
\rondu^{-1}$ the coarsest uniformity finer than $\rondu$ and
$\rondu^{-1}.$ The quasi-uniform space $(X,\rondu)$ is called {\it
bicomplete} if $(X,\rondu_s)$ is a complete uniform space. This
notion is useful and easy to handle, because one can appeal to
well known results  from the theory of uniform spaces.

A subset $Y$ of a quasi-uniform space $(X,\rondu)$ is called  {\it
precompact} if for every $U\in \rondu$ there  exists a finite
subset $Z$ of $X$ such that $Y\subset U[Z].$ The set $Y$ is called
{\it totally bounded } if  for every $U$ there exists a finite
family $A_1,...,A_n$ of subsets of $X$ such that $A_i\times
A_i\subset U,\, i=1,...,n,\,$ and  $Y\subset \cup_{i=1}^n A_i.$ In
uniform spaces total boundedness and precompactness agree, and a
set is compact if and only if it is totally bounded and complete.
A subset $Y$ of  quasi-uniform space $(X,\rondu)$ is totally
bounded if and only if it is totally bounded as a subset of the
uniform space $(X,\rondu_s).$

Another notion of completeness is that considered  by Sieber and
Pervin \cite{sieb-perv65}. A filter $\mathcal F$ in a
quasi-uniform space $(X,\rondu)$ is called $\rondu$-{\it Cauchy}
if for every $U\in \rondu$ there exists $x\in X$ such that
$U(x)\in \mathcal F.$ In terms of nets, a net $(x_\alpha,\alpha\in
D)$  is called $\rondu$-{\it Cauchy} if for every $U\in \rondu$
there exists $x\in X$ and $\alpha_0\in D$ such that
$(x,x_\alpha)\in U$ for all $\alpha\geq \alpha_0.$ The
quasi-uniform space $(X,\rondu)$ is called $\rondu$-{\it complete}
if every $\rondu$-Cauchy filter (equivalently, every
$\rondu$-Cauchy net) has a cluster point. If every such filter
(net) is convergent, then the quasi-uniform space $(X,\rondu)$ is
called $\rondu$-{\it convergence complete}. Obviously that
convergence complete implies complete, but the converse is not
true. It is clear that this notion corresponds to that of
$\rho$-completeness of a quasi-metric space. It is worth to notify
that the $\rondu_\rho$-completeness of the associated
quasi-uniform space $(X,\rondu_\rho)$ implies the
$\rho$-sequential completeness of the quasi-metric space
$(X,\rho),$ but the converse is not true (see \cite{kunz-reil93}).
The equivalence holds for the notion of left $K$-completeness
(which  will be defined immediately): a quasi-metric space is left
$K$-sequentially complete if and only if its induced
quasi-uniformity $\rondu_\rho$ is left $K$-complete
(\cite{rom92}).

A filter $\mathcal F$ in a
 quasi-uniform space $(X,\rondu)$ is
called {\it left $K$-Cauchy} provided for every $U\in \rondu$
there exists $F\in \mathcal F$ such that $U(x)\in F$ for all $x\in
F$. A net $(x_\alpha,\alpha\in D)$ in $X$ is called {\it left
$K$-Cauchy} provided for every $U\in \rondu$ there exists
$\alpha_0\in D$ such that $(x_\alpha,x_\beta)\in U$ for all
$\beta\geq \alpha\geq \alpha_0.$ The quasi-uniform space
$(X,\rondu)$ is called {\it left $ K$-complete} if every left
$K$-Cauchy filter (equivalently, every left $K$-Cauchy net)
converges.  If every left $K$-Cauchy filter converges with respect
to the uniformity $\rondu_s,$  then the quasi-uniform space
$(X,\rondu)$ is called {\it Smyth complete} (see \cite{kunz93b}
and \cite{smyth94}). This notion of completeness has applications
to computer science, see \cite{shek95}. In fact, there are a lot
of applications of quasi-metric spaces, asymmetric normed spaces
and quasi-uniform spaces to computer science, abstract languages,
complexity, see, for instance,
\cite{rom-raf-per02b,rom-raf-per04a,shek04,rom-shek99,rom-shek00b,rom-shek02}.

Another useful notion of completeness was considered  by
Doitchinov \cite{doich86,doich88,doich91a,doich91b,doich92}.  A
filter $\mathcal F$ in a quasi-uniform space $(X,\rondu)$ is
called $D$-{\it Cauchy} provided there exists a co-filter
$\mathcal G$ in $X$ such that for every $U\in \rondu$ there are
$G\in \mathcal G$ and $F\in \mathcal F$ such that $F\times
G\subset U. $ The quasi-uniform space $(X,\rondu)$ is called
$D$-{\it complete} provided every $D$-Cauchy filter converges. A
related notion of completeness was  considered by Andrikopoulos
\cite{andrik04a}. For a comparative study of the completeness
notions defined by pairs of filters see \cite{deak91} and
\cite{andrik04b}.

Notice also that these notions of completeness can be considered
within the framework of bitopological spaces in the sense of Kelly
\cite{kelly63}, since a quasi-metric space is a bitopological
space with respect to the topologies $\tau(\rho)$ and $\tau(\bar
\rho).$ For this approach see the papers by Deak
\cite{deak95,deak96}. It seems that  $K$ in the definition of left
$K$-completeness comes from Kelly who considered first this notion
(see \cite{zimm05}).

\section{Compact operators}

Recall that a subset $Z$ of an  asymmetric normed  space $(X,p)$
is called $p$-{\it precompact} if for every $\epsilon > 0$ there
exist
$z_1,...,z_n\in Z$ such that %
\bequ\label{def.t-bd} %
\forall z\in Z,\; \exists i\in \{1,...,n\},\; \;\;\; p(z-z_i)\leq \epsilon, %
\eequ %
or, equivalently, %
$$%
Z\subset U_\epsilon[\{z_1,...,z_n\}], %
$$%
where $U_\epsilon$ is the entourage %
$$%
U_\epsilon =\{(x,x')\in X\times X : p(x'-x)\leq \epsilon\} %
$$%
in the quasi-uniformity $\rondu_p\,.$

One obtains an equivalent notion taking the points $z_i$ in $X$
or/and $<\epsilon$ in \eqref{def.t-bd}.

 Let $(X,p),(Y,q)$ be asymmetric normed spaces and, as
before, let%
\bequ\label{def.rs} %
\mu \in \{p,\bar p,p_s\}\;\;\mbox{and}\;\; \nu \in \{q,\bar q,q_s\}. %
\eequ %
 A linear operator $A:X\to Y$ is called $(\mu,\nu)$-{\it compact} if the
set $A(B_\mu)$ is $\nu$-precompact in $Y.$

Some properties of compact operators are collected in the
following proposition. We shall denote by $(X,Y)^k_{\mu,\nu}$ the
set of all linear $(\mu,\nu)$-compact operators from $X$ to $Y.$
Notice that, for $\,\mu=p_s\,$ and $\,\nu=q_s,\,$ the space
$(X,Y)^\flat_{p_s,q_s}$ agrees with $(X,Y)^*_s,$ the
$(p_s,q_s)$-compact operators are the usual linear compact
operators between the normed spaces $(X,p_s)$ and $(Y,q_s),$ so
the proposition contains some  well
known results for compact operators on normed spaces.%
\begin{prop}\label{p.k-op1} %
Let  $(X,p),(Y,q)$ be asymmetric normed spaces. The following
assertions hold. %
\begin{enumerate} %
\item $(X,Y)^k_{\mu,\nu}$ is a semilinear subspace of $\,(X,Y)^\flat_{\mu,\nu}. $
\item $(X,Y)^k_{p,q}$ is $\tau(p,\bar q)$-closed in
$\,(X,Y)^\flat_{p,q}. $
\end{enumerate} %
\end{prop} %
\begin{proof} %
(1)\; We give the proof in the case $\mu=p$ and $\nu=q.$ The other
cases can be treated similarly.

If  $A:X\to Y$ is $(p,q)$-compact, then there exists
$\,x_1,...,x_n\in B_p\,$ such that %
\bequ\label{eq1.k-op1}%
\forall x\in B_p,\; \exists i\in \{1,...,n\},\;\quad q(Ax-Ax_i)\leq 1. %
\eequ%
If for $x\in B_p,\; i\in \{1,...,n\}\,$    is chosen
according to \eqref{eq1.k-op1}, then %
$$%
q(Ax)\leq q(Ax-Ax_i)+q(Ax_i)\leq 1+\max\{q(Ax_j): 1\leq j\leq n\}, %
$$%
showing that  the operator $A$ is $(p,q)$-bounded.

Suppose that $A_1,A_2:X\to Y$ are $(p,q)$-compact and let
$\epsilon >0.$  By the $ (p,q)$-compactness of the operators
$A_1,A_2,$ there exist $x_1,...,x_m$ and $y_1,...,y_n$ in $B_p$
such that %
$$%
\forall x\in B_p,\; \exists i\in \{1,...,m\},\; \exists j\in
\{1,...,n\},\;\; q(A_1x-A_1x_i)\leq \epsilon
\;\;\mbox{and}\;\;q(A_2x-A_2x_j)\leq \epsilon. %
$$%

It follows that for every $x\in B_p$ there exists a pair  $(i,j)$
with $1\leq i\leq m$ and $1\leq j\leq n$ such that %
$$%
q(A_1x+A_2x -A_1x_i-A_2y_j)\leq q(A_1x -A_1x_i)+q(A_2x
-A_2y_j)\leq 2\epsilon, %
$$%
showing that $\{A_1x_i+A_2y_j : 1\leq i\leq m,\; 1\leq j\leq n\}$
is a finite $2\epsilon$-net for $(A_1+A_2)(B_p).$

The proof of the compactness  of $\alpha A,$ for $\alpha > 0$ and
$A$ compact, is immediate and we omit it.\\

(2) {\it The $\tau(p,\bar q)$-closedness of}
$\,(X,Y)^k_{p.q}\,.$\\

Let $(A_n)$ be a sequence in  $(X,Y)^k_{p,q}$ which is
$\tau(p,\bar q)$-convergent to $A\in (X,Y)^\flat_{p,q}\,.$

For $\epsilon > 0$ choose $n_0\in \Nat$ such that %
\bequ\label{eq2.k-op1} %
\forall n\geq n_0,\;\forall x\in B_p,\;\; \bar q(A_nx-Ax)\leq
\epsilon\;\; (\iff q(Ax-A_nx)\leq\epsilon) . %
\eequ %

Let $x_1,...,x_m\in B_p$ be such that the points
$A_{n_0}x_i,\;1\leq i\leq m, \,$ form  an $\epsilon$-net for
$A_{n_0}(B_p).$  Then for every
$x\in B_p$ there exists $i\in \{1,...,m\}$ such that %
$$%
q(A_{n_0}x-A_{n_0}x_i)\leq \epsilon, %
$$%
so that, by \eqref{eq2.k-op1}, %
$$%
q(Ax-Ax_i)\leq q(Ax-A_{n_0}x)+q(A_{n_0}x-A_{n_0}x_i)+
q(A_{n_0}x_i-Ax_i) \leq 3\epsilon. %
$$%
Consequently, $Ax_i, \; 1\leq i\leq m,\,$ is a $3\epsilon$-net
for $A(B_p),$ showing that $A\in (X,Y)^k_{p,q}\,. $ %
\end{proof} %
\begin{rem} %
{\rm The assertion  (2) of Proposition \ref{p.k-op1} holds for
other types of compactness too, i.e. for the spaces
$(X,Y)^k_{\mu,\nu}$ with $\mu,\nu$ as in  \eqref{def.rs}, with
similar proofs.}
\end{rem} %

\section{The dual of a bounded linear operator}

Let $(X,p),\,(Y,q)$ be asymmetric normed spaces and $\mu,\nu$ as
in \eqref{def.rs}. For $A\in (X,Y)^\flat_{\mu,\nu}$ define
$A^\flat:Y^\flat_\nu\to X^\flat_\mu$  by%
\bequ\label{def.flatA}  %
A^\flat\psi =\psi\circ A,\;\; \psi\in Y^\flat_\nu. %
\eequ %

Obviously that $A^\flat$ is properly defined, additive and
positively homogeneous. Concerning the continuity we have. %
\begin{prop}\label{p.flatA1} %
\begin{enumerate} %
\item The operator $A^\flat$ is quasi-uniformly continuous with respect to
the quasi-uniformities  $\rondu_\nu^\flat$ and $\rondu_\mu^\flat$
on $Y^\flat_\nu$ and $X^\flat_\mu\,,$ respectively.
\item The operator $A^\flat$ is also
quasi-uniformly continuous with respect to the
$w^\flat$-quasi-unifor\-mities on $Y^\flat_\nu$ and
$X^\flat_\mu\,.$
\end{enumerate}
\end{prop} %
\begin{proof} %
(1)\; Take again $\mu=p$ and $\nu=q.$ For $\epsilon >0$ let %
$$%
U_\epsilon=\{ (\vphi_1,\vphi_2)\in X^\flat_p\times X^\flat_p :
\vphi_2(x)-\vphi_1(x)\leq \epsilon,\; \forall x\in B_p\}. %
$$%
 If $\|A|_{p,q}=0,$ then $A=0, $ so we can suppose
 $\|A|=\|A|_{p,q}>0.$  Let %
 $$%
V_\epsilon=\{ (\psi_1,\psi_2)\in Y^\flat_q\times Y^\flat_q :
\psi_2(x)-\psi_1(x)\leq \epsilon/\|A|,\; \forall x\in B_q\}. %
$$%

Taking into account that %
$$%
\forall x\in B_p,\; \vphi_2(x)-\vphi_1(x)\leq \epsilon/r\;\iff\;
\forall x'\in rB_p,\; \vphi_2(x')-\vphi_1(x')\leq \epsilon, %
$$%
and %
$$%
\forall x\in B_p,\;\quad q(Ax)\leq \|A|p(x) \leq \|A|, %
$$%
it follows %
$$%
A^\flat\psi_2(x)-A^\flat\psi_1(x)=\psi_2(Ax)-\psi_1(Ax)\leq
\epsilon, %
$$%
for all $x\in B_p, $ proving the quasi-uniform continuity of $A.$

(2) \; For   $x_1,...,x_n\in X$ and $\epsilon >0$   let %
$$%
V=\{(\vphi_1,\vphi_2)\in X^\flat_p\times X^\flat_p :
\vphi_2(x_i)-\vphi_1(x_i)\leq \epsilon,\; i=1,..., n\} %
$$%
be a $w^\flat$-entourage in $X^\flat_p.$ Then %
$$%
U=\{(\psi_1,\psi_2)\in Y^\flat_q\times Y^\flat_q :
\psi_2(Ax_i)-\psi_1(Ax_i)\leq \epsilon,\; i=1,...,n\}, %
$$%
is a $w^\flat$-entourage in $  Y^\flat_q$ and %
$(A^\flat\psi_1,A^\flat\psi_2)\in V$ for every $(\psi_1,\psi_2)\in
U, $ proving the quasi-uniform continuity of $A^\flat$ with
respect to the $w^\flat$-quasi-uniformities on $Y_q^\flat$ and
$X^\flat_p\,$\,. %
\end{proof} %

Now we can prove the analog of the Schauder theorem for the
asymmetric dual. %
\begin{theo}\label{t.Sch} %
Let $(X,p), (Y,q)$ be  asymmetric normed spaces. If the linear
operator $A:X\to Y$ is $(p,q)$-compact, then $A^\flat(B^\flat_q)$
is precompact with respect to the quasi-uniformity $\rondu_p^\flat$ on $X^\flat_p.$ %
\end{theo} %
\begin{proof} %
For  $\epsilon > 0 $ let %
$$%
U_\epsilon=\{(\vphi_1,\vphi_2)\in X^\flat_p\times X^\flat_p :
\vphi_2(x)-\vphi_1(x)\leq \epsilon, \; \forall x\in B_p\}, %
$$%
be an entourage in $X^\flat_p$ for the quasi-uniformity
$\rondu_p^\flat\,.$

Since $A$ is $(p,q)$-compact, there exist
$x_1,...,x_n \in B_p\,$ such that %
\bequ\label{eq1.Sch} %
\forall x\in B_p\,,\; \exists i\in \{1,...,n\},\;\quad
q(Ax-Ax_i)\leq \epsilon. %
\eequ %

By the Alaoglu-Bourbaki theorem (\cite[Theorem
4]{rom-raf-per03a}), the set $B^\flat_q$ is $w^\flat$-compact, so
by the $(w^\flat,w^\flat)$-continuity of the operator $A^\flat$
(Proposition \ref{p.flatA1}),  the set $A^\flat(B^\flat_q)$ is
$w^\flat$-compact in $X^\flat_p.$ Consequently, the $w^\flat$-open
cover of $A^\flat(B^\flat_q),$%
$$%
V_\psi=\{\vphi\in X^\flat_p : \vphi(x_i)-A^\flat\psi(x_i)
<\epsilon,\; i=1,...,n\},\;\psi\in B^\flat_q, %
$$%
contains a finite subcover $V_{\psi_k},\, 1\leq k\leq m, $ i.e,
\bequ\label{eq2.Sch} %
A^\flat(B^\flat_q)\subset \bigcup\{V_{\psi_k} : 1\leq k\leq m\}, %
\eequ %
for some $m\in \mathbb N$ and $\psi_k\in B^\flat_q,\; 1\leq k\leq
m.$

Now let  $\psi \in B^\flat_q\,.$  By \eqref{eq2.Sch} there exists
$k\in \{1,...,m\}$ such that %
$$%
A^\flat\psi(x_i)-A^\flat\psi_k(x_i) < \epsilon,\; i=1,...,n. %
$$%
 If  $x\in B_p,$ then, by \eqref{eq1.Sch}, there  exists $i\in
 \{1,...,n\}, $ such that %
 $$%
 q(Ax-Ax_i)\leq \epsilon. %
 $$%

 It follows  %
 $$%
 \begin{aligned} %
 &\psi(Ax)-\psi_k(Ax) =\\
 & =
 \psi(Ax)-\psi(Ax_i)+\psi(Ax_i)-\psi_k(Ax_i)+\psi_k(Ax_i)-\psi(Ax_i)\\
&\leq 2q(Ax-Ax_i)+\epsilon \leq 3\epsilon. %
\end{aligned} %
$$%

Consequently, %
$$%
\forall x\in B_p\,,\; \quad (A^\flat\psi-A^\flat\psi_k)(x)\leq 3\epsilon, %
$$%
proving that %
$$%
A^\flat(B^\flat_q)\subset U_{3\epsilon}[\{A^\flat\psi_1,...,A^\flat\psi_m\}]. %
$$%
\end{proof} %

{\bf Comments}\; As a  precaution, we have defined the compactness
of an operator $A$ in terms of the precompactness of the image of
the unit ball $B_p \,$ by $A,$ rather than by the relative
compactness of $A(B_p)\,,$ as in the case of compact operators on
usual normed spaces. As can be seen from Section 2, the relations
between precompactness, total boundedness and completeness are
considerably more complicated in the asymmetric case than in the
symmetric one. To obtain some compactness properties  of the set
$A(B_p)\,,$ one needs a study of the completeness of the space
$(X,Y)^\flat_{\mu,\nu}$ with respect to various quasi-uniformities
and various notions of completeness, which could be the topic of
further investigation.

\providecommand{\bysame}{\leavevmode\hbox
to3em{\hrulefill}\thinspace}
\providecommand{\MR}{\relax\ifhmode\unskip\space\fi MR }
\providecommand{\MRhref}[2]{%
  \href{http://www.ams.org/mathscinet-getitem?mr=#1}{#2}
} \providecommand{\href}[2]{#2}

\end{document}